\documentclass[11pt,a4paper]{amsart}

\usepackage{amsmath,amsthm}

\theoremstyle{plain}
\newtheorem{thm}{Theorem}[section]
\newtheorem{cor}[thm]{Corollary}
\newtheorem{prop}[thm]{Proposition}
\newtheorem{lem}[thm]{Lemma}
\newtheorem{ques}[thm]{Question}

\theoremstyle{definition}
\newtheorem{defn}[thm]{Definition}
\newtheorem{exmp}[thm]{Example}

\theoremstyle{remark}
\newtheorem{rem}{Remark}[section]

\let\lvert=|\let\rvert=|

\usepackage{amsfonts}

\newcommand{\calL}{{\mathcal L}}

\newcommand{\calA}{{\mathcal A}}

\newcommand{\e}{\epsilon}
\newcommand{\clco}{\overline{co}}
\newcommand{\wclco}{\overline{co}^{w^{\ast}}}

\newcommand{\ospan}{\overline{\mbox{span}}}
\newcommand{\xast}{x^{\ast}}

\newcommand{\yast}{y^{\ast}}
\newcommand{\Yast}{Y^{\ast}}

\newcommand{\Xast}{X^{\ast}}

\newcommand{\Xastast}{X^{\ast\ast}}

\newcommand{\xastast}{x^{\ast\ast}}

\newcommand{\Tast}{T^{\ast}}
\newcommand{\Sast}{S^{\ast}}

\newcounter{abcd}
\newcounter{iiiii}
\newenvironment{ekvivalens}
{\setcounter{iiiii}{0}
\begin{list}%
{{\rm (\roman{iiiii})}} 
{\usecounter{iiiii}
\parsep=\parskip
\topsep=1pt plus 2pt minus 1pt
\itemsep=1pt plus 2pt minus 1pt
\leftmargin=3\baselineskip
\labelsep=.6\baselineskip
\labelwidth=2.4\baselineskip
\rightmargin 0pt}%
}%
{\end{list}}
\newenvironment{statements}%
{\setcounter{abcd}{0}
\begin{list}%
{{\rm (\alph{abcd})}} 
{\usecounter{abcd}
\parsep=\parskip
\topsep=1pt plus 2pt minus 1pt
\itemsep=1pt plus 2pt minus 1pt
\leftmargin=3\baselineskip
\labelsep=.6\baselineskip
\labelwidth=2.4\baselineskip
\rightmargin 0pt}%
}%
{\end{list}}
\begin{document}
\title[Boundedness and surjectivity]{Boundedness and surjectivity in normed spaces}
\author{Olav Nygaard}

\address{Deparment of Mathematics, Agder College, Tordenskjoldsgate 65, 4604 Kristiansand, Norway }
\email{Olav.Nygaard@hia.no}
\thanks{Filename: massive.tex. AmsLaTeX.}
\date{\today}

\begin{abstract} We define the ($w^\ast$-) boundedness property and the ($w^\ast$-) surjectivity property for sets in normed spaces. We show that these properties are pairwise equivalent in complete normed spaces by characterizing them in terms of a category-like property called ($w^\ast$-)  thickness. We give examples of interesting sets having or not having these properties. In particular, we prove that the tensor product of two $w^\ast$-thick sets in $\Xastast$ and $\Yast$ is a $w^\ast$-thick subset in $L(X,Y)^\ast$ and obtain as a concequense that the set $w^\ast -exp\:B_{K(l_2)^\ast}$ is $w^\ast$-thick.
\end{abstract}

\maketitle

\vspace*{5mm}
\centerline{\sc Introduction}

Recall the Banach-Steinhaus theorem for Banach spaces: \em A family of linear continuous operators on a Banach space $X$, which is pointwise bounded on a set of second category, is bounded. \em

Let $X$ be a normed linear space. Motivated by the Banach-Steinhaus theorem we say that $A\subset X$ has the \em boundedness property \em  if every family of linear continuous operators on $X$, which is pointwise bounded on $A$, is bounded. More generally, if $Y$ is a normed space and $\calA$ is a subset of $\calL(X,Y)$, we say that $A$ has the \em $\calA$- restricted boundedness property \em  if every family of linear continuous operators in $\calA$, which is pointwise bounded on $A$, is bounded. In the latter definition, if $\calA$ is the space of adjoints, we say that $A\subset \Xast$ has the \em $w^\ast$-boundedness property \em .  

From the proof of the Banach-Steinhaus theorem we conclude that every second category set $A$ in a Banach space $X$ has the boundedness property. However, the Nikodym-Grothendieck boundedness theorem (see e.g. \cite[p. 14]{DU} or \cite[p. 80]{Diestel}) says in our terminology exactly that the set of characteristic functions on the unit sphere of $B(\Sigma)$ has the boundedness property. This set is certainly not of the second category, it is even nowhere dense. Thus it may be possible to sharpen the Banach-Steinhaus theorem.

Let us have a look at a more recent theorem of J. Fernandez \cite{fernandez} (see also \cite{shapiro}), which is in the same spirit as the Nikodym-Grothendieck theorem.
\newpage

\begin{thm}\label{fernandez} Suppose $(f_n)$ is a sequence in $L_1 (T)$ such that 
\[|\sup_n \int_T f_n \phi\:d\theta|<\infty\]
for all inner functions $\phi$. Then 
\[|\sup_n \int_T f_n g\:d\theta|<\infty\]
for all $g\in H^\infty (T)$ (and hence $(f_n)$ is bounded).
\end{thm}

It is well known that the pre-dual of $H^\infty$ is $L_1 /H_0^1$. Thus, by our Theorem \ref{wunicar}, in our language this theorem says that the set of inner functions has the $w^\ast$-boundedness property in $H^\infty (T)$. In \cite{fernandez} and \cite{fhs} the question whether the set of inner functions has the boundedness property was posed. In \cite{hui} it was shown that also the set of Blaschke-products has the $w^\ast$-boundedness property in $H^\infty (T)$.

Now recall one part of the Open mapping theorem: \em Suppose $Y$ is a Banach space and $T\in\calL (Y,X)$. If $A\subset X$ is of the second category in $X$ and $TY\supset A$, then $TY=X$. \em 

Motivated by this we say that a set $A\subset X$ has the \em surjectivity property \em if, for every Banach space $Y$, every $T:Y\rightarrow X$ onto $A$ is onto $X$. If the conclusion holds for a subset $\calA \subset\calL (X,Y)$, we say that $A$ has the \em $\calA$-restricted surjectivity property \em . A special case of this is the case where $\calA$ is the space of adjoint operators between two duals. In this case we will say that $A$ has the \em $w^\ast$-surjectivity property \em . 

Again second category sets in Banach spaces have the surjectivity property. And again the set of characteristic functions on the unit sphere of $B(\Sigma)$ has the
surjectivity property. This follows directly from a theorem of Seever (see \cite{seever} or \cite[p. 17]{DU}). 

A natural question arises: Will the set of inner functions in $H^\infty (T)$ have the $w^\ast$- surjectivity property? We will see that the ($w^\ast$-) surjectivity property and the ($w^\ast$-) boundedness property are equivalent in Banach spaces and hence the answer to this question is "yes". This also shows that the Nikodym-Grothendieck theorem and Seever's theorem are closely related. The Nikodym-Grothendieck theorem has a well of consequences and one might expect nice consequences of analog theorems in other Banach spaces as well. In the following sections we will show different techniques to obtain such theorems. 

\section{Thick and thin sets}

One objective in this paper is to prove that in Banach spaces the ($w^\ast$-) surjectivity property and the ($w^\ast$-) boundedness property are equivalent to a common property, called ($w^\ast$-) thickness.

For this concepts we now make a list of five fundamental properties of sets in normed spaces:
\begin{defn}\label{egenskaper}
\begin{ekvivalens}
\item A set $A\subset X$ such that $\ospan\:A =X$ is called \em fundamental\em.
\item A set $A\subset X$ such that $\inf_{f\in S_{\Xast}}\sup_{x\in A}|f(x)|\geq\delta$, for some $\delta >0$, is called \em norming (for $\Xast$)\em.
\item A set $A\subset X$ such that $tA\supset B_X$ for some $0<t<\infty$ is called \em absorbing \em.
\item A set $A\subset X$ such that for all $\e>0$, there exists $0<t<\infty$ such that $tA+\e B_X \supset B_X$ is called \em almost absorbing\em.
\item A set $A\subset X$ such that there exist $0<\lambda <1$ and $0<t<\infty$ such that $tA +\lambda A \supset B_X$ is called \em $\lambda$-almost absorbing\em.  
\end{ekvivalens}
\end{defn}

By the Hahn-Banach separation theorem the following lemma is easy to prove (see Remark \ref{complex} for the complex case).
\begin{lem}\label{norm1} The following statements are equivalent.
\begin{statements}
\item $A$ is norming for $\Xast$.
\item $\clco(\pm A)$ is norming for $\Xast$.
\item There exists a $\delta>0$ such that $\clco(\pm A)\supseteq\delta B_X$, i.e. $\clco(\pm A)$ is absorbing.
\end{statements}
\end{lem}

More specifically, we can speak about $\delta$-norming sets, where the $\delta$ refers to the $\delta$ in (c).

Suppose a set $B\in \Xast$ is such that $\inf_{x\in S_{X}}\sup_{f\in B}|f(x)|\geq\delta$ for some $\delta >0$. In this case we will call the set \em norming for $X$ \em or \em $w^\ast$-norming \em . Of course we have a similar lemma for sets which are norming for $X$. 
\begin{lem}\label{norm2} The following statements are equivalent.
\begin{statements}
\item $B$ is norming for $X$.
\item $\clco(\pm B)$ is norming for $X$.
\item $\wclco(\pm B)$ is norming for $X$.
\item There exists a $\delta>0$ such that $\wclco(\pm B)\supseteq\delta B_{\Xast}$, i.e. $\wclco(\pm B)$ is absorbing.
\end{statements}
\end{lem}

\begin{prop} In normed spaces exactly the following implications are valid. (a) $\Rightarrow$ (b) $\Rightarrow$ (c) $\Rightarrow$ (d) $\Rightarrow$ (e).
\begin{statements}
\item $A$ is absorbing.
\item $A$ is almost absorbing.
\item $A$ is $\lambda$-almost absorbing.
\item $A$ is norming.
\item $A$ is fundamental.
\end{statements}
For convex sets (c) $\Leftrightarrow$ (b). For closed, convex, symmetric sets (a) $\Leftrightarrow$ (d). But (e) does not imply (d) even in this case.
\end{prop}

\begin{proof} (a) $\Rightarrow$ (b) $\Rightarrow$ (c) is trivial. If $A$ is $\lambda$-absorbing, there exists $\lambda, t$ such that $B_X\subset tA+\lambda B_X$. Then, of course $B_X\subset t\cdot\clco(\pm A) + \lambda B_X$. Thus 
\[B_X\subset t\cdot\clco(\pm A) + \lambda B_X \subset t\cdot\clco(\pm A) + \lambda    (t\cdot\clco(\pm A) + \lambda B_X)\subset\cdots\]
\[\subset t\cdot\clco(\pm A)+\lambda t\cdot\clco(\pm A)+\cdots +\lambda^n t\cdot\clco(\pm A)+\lambda^{n+1} B_X.\]
Since $\clco(\pm A)$ is convex we obtain
\[B_X\subset t\cdot\clco(\pm A)(1+\lambda+\cdots\lambda^n )+\lambda^{n+1} B_X.\]
But this is true for every $n$. Hence we can take limits:
\[B_X \subset t\cdot\clco(\pm A) \frac{1}{1-\lambda},\]
which gives
\[\clco(\pm A)\supset\frac{1-\lambda}{t}B_X,\]
and by Lemma \ref{norm1}, $A$ is norming. That (d) $\Rightarrow$ (e) is trivial.

To see that (d) does not imply (c) and (c) does not imply (b) take $A=\{e_n\}$, the unit vectors in $l_1$. 

We now prove that (c) implies (b) for convex sets. Let $\e>0$. Since $A$ is $\lambda$-almost absorbing and convex
\[B_X \subset tA+\lambda B_X \subset tA+\lambda(tA+\lambda B_X)=t(1+\lambda)A+\lambda^2 B_X\]
\[\subset \cdots \subset (1+\lambda +\lambda^2+ \cdots +\lambda^{n-1})tA +\lambda^{n}B_X,\]
for every natural number $n$. Since $\lambda<1$, $\lambda^n$ is eventually less than $\e$.

To see that (b) does not imply (a) take a dense subset of the unit ball in a separable Banach space. That (a) $\Leftrightarrow$ (b) $\Leftrightarrow$ (c) $\Leftrightarrow$ (d) for closed convex and symmetric sets follows from Lemma \ref{norm1}.

To end the proof, by Lemma \ref{norm1}, it is enough to construct a closed, convex, symmetric set $A$ such that $A$ is fundamental but not norming. For this, let 
\[A=\clco(\pm \frac{e_n}{n})\subset l_1.\]
It is easy to check that $A$ has the desired properties.
\end{proof}

Later we will need the now well-known construction of Davis, Figiel, Johnson and Pe{\l}zcynski. See \cite{DFJP} or \cite[p. 227-228]{Diestel}, or \cite{nylioja} for some recent results. We call the space constructed by this procedure from a bounded, absolutely convex set $K$, the DFJP-space constructed on $K$.

\begin{prop}\label{ektedfjp} Let $K\subset X$ be bounded, convex and symmetric. Let $X_K$ be the DFJP-space constructed on $K$. Then the embedding $J_K$ of $X_K$ into $X$ is an isomorphism if and only if $K$ is norming for $\Xast$. 
\end{prop} 

\begin{proof} If $K$ is norming for $X$, $\overline{J_K (K)}\supset\lambda\cdot B_X$ for some $\lambda>0$. Thus $\overline{J_K (B_{X_K})}\supset\lambda\cdot B_X$ and $J_K$ is invertible. 

If $X_K$ and $X$ are isomorphic by $J_K$ then, for some $\delta>0$,
\[\delta B_X\subset C \stackrel{\rm def}{=} J_K (B_{X_K}).\]
But, by the construction of $X_K$,
\[C\subset a^n K + a^{-n}B_X\]
for all $n$. We use that $\delta B_X\subset C$ inductively for a constant $n$:
\[a^n K + a^{-n}B_X\subset a^n K+\frac{a^{-n}}{\delta}(a^n K+a^{-n}B_X)\]
\[\subset (a^n +\frac{1}{\delta})K+\frac{a^{-2n}}{\delta}B_X
\subset(a^n +\frac{1}{\delta})K+\frac{a^{-2n}}{\delta^2}(a^n K+a^{-n}B_X)\]
\[\subset(a^n +\frac{1}{\delta}+\frac{a^{-n}}{\delta^2})K+\frac{a^{-3n}}{\delta^2}B_X
\] 
Continuing this way gives after $r$ steps
\[C\subset\left(a^n +a^n \sum_{k=1}^{r}\frac{1}{(\delta\cdot a^n )^k}\right)K + \frac{a^{-(r+1)n}}{\delta^r}\]
\[=\left(a^n +a^n \sum_{k=1}^{r}\frac{1}{(\delta\cdot a^n )^k}\right)K + (\frac{1}{\delta a^n})^r \cdot a^{-n}\]
This is true for any $n$. Now choose $n$ so big that $\delta a^n >2$. Then, by letting $r\rightarrow\infty$, we obtain
\[\frac{\delta}{2a^n} B_X \subset \overline{K}.\]
This proves that $K$ is norming.
\end{proof}

\begin{rem}\label{complex} When the space under consideration is complex the definition of a norming set is of course that \em for some $\delta>0$, $\clco(rA)\supset\delta B_X$ for all $|r|=1$. \em It is easy to verify that all the results so far are true with complex scalars instead of reals. Note also that when testing for the boundedness property and the surjectivity property one could as well test on the set $co(rA),\:|r|=1$. Thus, for the rest of the paper we can assume $A$ to be symmetric and convex.
\end{rem}

Fundamental sets are useful when testing for weak-star and weak convergence of nets. The following Proposition is well-known and classic, but we feel that it deserves some space here. The proof is of course elementary.

\begin{prop} Let $X$ be a normed space. 
\begin{statements}
\item Suppose a \em bounded \em net in $\Xast$ converges pointwise on a fundamental set $A\subset X$. Then it converges weak-star.
\item Suppose a \em bounded \em net in $X$ converges pointwise on a fundamental set $A\subset \Xast$. Then it converges weakly.
\end{statements}
\end{prop}

\begin{rem} The Simons-Rainwater theorem \cite[e.g.]{simons} is trivial (and works also for bounded nets) whenever the James boundary under consideration is fundamental in $\Xast$ (e.g. when $\Xast$ has the RNP).
\end{rem}

\begin{exmp} To illustrate the terms fundamental and norming a bit more, we can look at Blaschke products in $H^\infty(T)$. A famous theorem of Marshall  says: \em $H^\infty$ is the closed, linear span of the Blaschke products. \em In other words, the Blaschke products form a fundamental subset of $H^\infty$. Later it has been shown that the Blaschke products form a subset of $H^\infty$ which is 1-norming for the dual (see \cite[p. 195-197]{garnett}). The most recent theorem in this direction is, as far as I know, the result from \cite{oyma} saying that the interpolating Blaschke products form a set which is $10^{-7}$-norming for the dual. 
\end{exmp}

We now define the term thick set: 
\begin{defn} A set is called ($w^\ast$-) \em thick \em if it is not an increasing union of non-($w^\ast$-)norming sets. A set which is not ($w^\ast$-) thick is called \em ($w^\ast$-) thin \em .
\end{defn}

This classification of sets is not standard and the terms thick and thin sets are often used to describe properties of sets. Maybe it would be better to call thick sets \em Fonf sets \em since I think he is the one who first and best demonstrated the relevance and importance of the thick sets. Fonf however never uses the word thick, but in his works he always operates with thin sets. See e.g. \cite{KF}, \cite{Fwe} and \cite{Fs} for examples of earlier use and applications of these concepts. 

\begin{exmp}\label{hjorne}
To get an idea of these properties, one can think of the extreme points $A$ of the unit ball in $l_1$. This is a countable, hence thin (and thus $w^\ast$-thin) set. $B_{l_1}$ is the norm-closure of the convex hull of its extreme points $A=ext\:B_{l_{1}}=\{\pm e_i\}_{i=1}^{\infty}$ so it is 1-norming. Let $f_n = ne_n $. Then $\{f_n\}\subset c_0$ is pointwise bounded on $A$, but obviously not bounded on all of $l_1$. So this $A$ does not have the $w^\ast$-boundedness property. Now define an operator $T\in L(l_\infty ,l_1)$ by $T(x_1 ,x_2, ...)=(2^{-n}x_n )_{n=1}^{\infty}$. Then $T$ is onto $A$, but since it is injective, it can't be onto $l_1$. Thus $A$ does not have the surjectivity property. Since $T$ is the adjoint of $S:c_0 \rightarrow l_1$ defined by $S(x_1 ,x_2, ...)=(2^{-n}x_n )_{n=1}^{\infty}$, $A$ does not have the $w^\ast$-surjectivity property either.
\end{exmp}

We want to show that the simple Example \ref{hjorne} is just a special case of very general principles of normed spaces.

\section{The boundedness property in normed spaces}

We recall the definition of the boundedness property.

\begin{defn} A subset $A$ of a normed linear space $X$ is said to have the boundedness property if for every normed space $Y$ every family $(T_\alpha)\subset\calL (X,Y)$, which is pointwise bounded on $A$, is bounded.
\end{defn}

A special variant of the following theorem was first published in \cite{nyg}. Also parts of it are implicit in \cite{Fs}.

\begin{thm}\label{unicar}Suppose $A$ is a subset of a normed space $X$. The following statements are equivalent.
\begin{statements}
\item $A$ has the boundedness property.
\item $A$ has the $\calA$-restricted boundedness property where $\calA$ is the collection of all countable sets of linear bounded functionals on $X$.
\item $A$ is thick.
\end{statements}
\end{thm}

\begin{proof} (a) clearly implies (b). Suppose $A$ is thin. Since $A$ is thin we can pick a countable, increasing covering, $\cup A_n$ of $A$, consisting of sets which are non-norming for $\Xast$, and a sequence $(f_n)\subset \Xast$ such that $f_n \in nS_{\Xast}$ but $\sup_{A_{n}}|f_n (x)|<1$. Let $x$ be an arbitrary element of $A$. Then there is a natural number $m$ such that $x\in B_m$. Thus, since $(A_n)$ is increasing,
\[|f_k (x)|\leq ||f_k||\:||x||<m||x||\hspace{1cm}\mbox{if}\: k < m,\]
while
\[|f_k (x)|\leq 1\hspace{1cm}\mbox{if}\: k\geq m.\]
Thus (b) implies (c).

To show that (c) implies (a), suppose $(T_{\alpha})$ is pointwise bounded on $A$, i.e.
\[\sup_{\alpha} \|T_{\alpha} x\|<\infty \hspace{1cm}\mbox{for all}\: x\in A.\]
Put $A_m = \{x\in A\::\sup_\alpha \|T_\alpha x\|\leq m\}.$ Then $(A_m)$ is an increasing family of sets which covers $A$. Since $A$ is thick, some $A_q$ is norming. Then, using Lemma \ref{norm1}, there exists a $\delta>0$ such that 
\[\overline{co}(\pm A_q)\supseteq \delta B_X.\]
But then, for arbitrary $\alpha$,
\[\delta \|T_\alpha\| = \sup_{x\in\delta S_X} \|T_\alpha x\|\leq\sup_{x\in\overline{co}(\pm A_q)} \|T_\alpha x\|\leq  q.\]
Thus $\sup_{\alpha}\|T_\alpha\|\leq q/\delta<\infty$ and the theorem is proved.
\end{proof}

\begin{rem} Note how simple and general Theorem \ref{unicar} is. It gives an easy and transparent way to the classical Banach-Steinhaus theorem for Banach spaces, as soon as we have shown that sets of second category are thick. This is done in the following trivial Lemma.
\end{rem}
 
\begin{lem} Let $X$ be a normed space and suppose $A$ is a second category set in $X$. Then $A$ is thick.
\end{lem}

\begin{proof} Suppose $A$ is covered by an increasing family $(A_i)$. Since $A$ is of second category some $\overline{A_m }$, contains a ball. Then $\clco(\pm A_m)$ contains a ball centered at the origin, and hence $A_m$ is norming. Since $(A_i)$ was arbitrary, $A$ must be thick.
\end{proof}

Of course we also have a "uniform boundedness theorem" characterizing $w^\ast$-thickness. It goes like this.

\begin{thm}\label{wunicar} Suppose $B$ is a bounded subset of a dual space $\Xast$. Then the following statements are equivalent.
\begin{statements}
\item $B$ has the $w^\ast$-boundedness property.
\item $B$ has the $\calA$-restricted boundedness property where $\calA$ is the collection of all sequences in $X$.
\item $B$ is $w^\ast$-thick.
\end{statements}
\end{thm}

\begin{proof} (a) implies of course (b). The proof that (b) implies (c) is completely analog to the corresponding part of the proof of Theorem \ref{unicar}. The proof that (c) implies (a) is also very similar to the corresponding part of the proof of Theorem \ref{unicar}. Just put
\[B_m =\{\xast\in B\::\sup_{\alpha}\|\Tast_{\alpha}\xast\|\leq m\}.\]
\end{proof}

\begin{rem} The set of extreme points of the unit ball of $l_1$ shows that a set can be norming for the dual without being $w^\ast$-thick. It is also possible for a set to be $w^\ast$-thick without being norming for the dual (although it is of course norming for the pre dual). In fact, the unit ball of any non-reflexive space, considered as a subset of the bidual, give examples of such situations (They are not even fundamental). The Blaschke products in $H^\infty$ is an example of a set which is both $w^\ast$-thick and norming for the dual space.
\end{rem}

\section{The surjectivity property in Banach spaces}

In this section we study the surjectivity property in Banach spaces. Let us start with the formal definition:

\begin{defn} In a normed linear space $X$ a set $A$ is said to have the surjectivity property if for every normed linear space $Y$, every $T\in\calL (Y,X)$, such that $TY\supset A$, is onto $X$. If $A$ in a normed space has the surjevtivity property for all $T$'s coming from Banach spaces, we say that $A$ has the surjectivity property for Banach spaces.
\end{defn} 

Recall that an operator $T:Y\rightarrow X$ is called Tauberian if $(T^{\ast\ast})^{-1}(X)\subset Y$. As an intuition, it is often helpful to think of these operators as opposite to weakly compact operators. A nice reference for the theory of Tauberian operators is \cite{gonz}.

Here is the connection between thickness and the surjectivity property. The theorem is an extension of a theorem discovered by Kadets and Fonf \cite{KF}.  

\begin{thm}\label{surcar} Suppose $A$ is a subset of a Banach space $X$. The following statements are equivalent.
\begin{statements}
\item $A$ has the surjectivity property for Banach spaces.
\item $A$ has the $\calA$- restricted surjectivity property where $\calA$ is the set of all injections from Banach spaces into $X$.
\item $A$ has the $\calA$- restricted surjectivity property where $\calA$ is the set of all Tauberian injections from Banach spaces into $X$.
\item $A$ is thick
\end{statements}
\end{thm}

\begin{proof}
Of course (a) implies (b) and (b) implies (c). 

To show that (c) implies (d) suppose (d) is not true, i.e. $A$ is thin. We will construct a Tauberian injection which is onto $A$ but not onto all of $X$. Let $(A_i)$ be an increasing family of subsets of $A$ such that $A=\cup_{i=1}^{\infty}A_i$. Since $\cup_{i=1}^{\infty}A_i
=\cup_{i=1}^{\infty}A_i \cap i\cdot B_X$, we may assume each $A_i$ to be contained in $i\cdot B_X$. Put $C_1=A_1$ and $C_i =A_i \setminus A_{i-1}$. Define
\[C=\clco(\pm\bigcup_{i=1}^{\infty}\frac{C_i}{i^{2}}).\]
Then $C$ is closed, bounded, convex and symmetric. We now show that $C$ is non-norming for $\Xast$. To do this, let $\e>0$ and take $j$ such that $1/j<\e$. Since $A_j$ is not a norming set, there is a functional $f\in S_{\Xast}$ such that $\sup_{x\in A_{j}}|f(x)|<\e$. But then, by the definition of $C$,
\[\sup_{x\in C}|f(x)|=\sup_i \left\{ \frac{1}{i^{2}}\sup_{x\in C_{i}}|f(x)|\right\}<\e.\]
Hence, by Proposition \ref{ektedfjp}, the Davis-Figiel-Johnson-Pelzcynski construction will produce a Banach space $Y$ and an operator $J:Y\rightarrow X$ with the desired properties, i.e. it is injective, Tauberian, onto $A$ but not onto all of $X$.

It remains to show that (d) implies (a). To do this, let $T$ be any bounded, linear operator into $X$ and onto $A$. Put $A_i = T(i\cdot B_Y)\cap A$, where $Y$ is the domain space of $T$. Since $T$ is onto $A$, $A=\cup_{i=1}^{\infty} A_i$, an increasing covering of $A$. Since $A$ is thick some $A_j$ is norming for $\Xast$. By Lemma \ref{norm1}, there exists a $\delta>0$ such that
\[j\cdot\clco(\pm TB_Y)=j\cdot\overline{TB_Y}\supset \delta B_X.\]
Hence $\overline{TB_Y} \supset (\delta/j)\cdot B_X$ and, by e.g. \cite[Thm 4.13]{rudin}, $T$ is onto. 
\end{proof}

\begin{cor}\label{equivcor} In Banach spaces, the surjectivity property for Banach spaces and the boundedness property are both equivalent to thickness.
\end{cor}

\begin{cor} $A\subset X$, where $X$ is a Banach space is thick if and only if the completion of any norming of \em span \em $A$, stronger than the original one, is $X$.
\end{cor}

\begin{rem}
Note that, since the characteristic functions is a thick set in $B(\Sigma)$, Seever's theorem follows as a special variant of the very general Theorem \ref{surcar}. 
\end{rem}

\begin{rem}
Also note that the DFJP-embedding $J$ is "a little more" than Tauberian. Also $J^{\ast\ast}$ is Tauberian. In \cite{nylioja} a very easy argument is given to show that, in fact, $J$ is a norm-norm homeomorphism when restricted to the set on which the DFJP-space is constructed.
\end{rem}

It is interesting to compare Theorem \ref{surcar} with the following observation of R. Neidinger \cite[p.119]{neiddoc}. It's proof is nothing but a close inspection of a standard proof of a general open mapping theorem, see e.g. \cite[p.48]{rudin}. An interesting application can be found in \cite{Neidinger}.

\begin{lem}\label{neidlem} Let $T\in\calL (X,Y)$. Then the following are equivalent.
\begin{statements}
\item $TB_X$ almost $\lambda$-absorbs $B_Y$.
\item $TB_X$ almost absorbs $B_Y$.
\item $TB_X$ absorbs $B_Y$.
\item $T$ is onto.
\end{statements}
\end{lem}

An analog to Lemma \ref{neidlem}, using the terms thick set and norming set is the following:

\begin{lem} Let $T\in\calL (X,Y)$. Then the following are equivalent.
\begin{statements}
\item $TB_X$ is norming for $\Yast$.
\item $TB_X$ is thick in $Y$.
\item $T$ is onto.
\end{statements}
\end{lem}

\begin{proof} Only the implication (a) $\Rightarrow$ (b) needs proof since the remaining implications are trivial. Suppose $TB_X$ is thin. Write $TB_X =\cup_{i=1}^\infty A_i,$
an increasing union of sets which are non-norming for $\Yast$. Then $B_X =\cup_{i=1}^\infty (T^{-1}(A_i)\cap B_X )$, an increasing union of sets. Since $B_X$ is thick, there exists a number $m$ and a $\delta >0$ such that $\clco(T^{-1}(A_m)\cap B_X )\supset\delta B_X$. Thus
\[\delta TB_X \subset T(\clco(T^{-1}(A_m)\cap B_X ))\subset \clco(T(T^{-1}(A_m)\cap B_X ))\]
\[\subset \clco(A_m \cap TB_X)=\clco(A_m).\]
Thus
\[\clco(TB_X)=\overline{TB_X}\subset\frac{1}{\delta}\clco(A_m).\]
This shows that $TB_X$ is non-norming.
\end{proof}

\begin{cor} Thin, norming sets are never linear images of unit balls.
\end{cor}

Let us now consider $w^\ast$-thick sets. Here is a characterization of such sets in terms of surjectivity properties:

\begin{thm}\label{wsurcar} Suppose $B$ is a subset of a dual space $\Xast$. Then the following statements are equivalent.
\begin{statements}
\item $B$ has the surjectivity property for all dual operators into $\Xast$, i.e. $B$ has the $w^\ast$- surjectivity property.
\item $B$ has the surjectivity property for all dual injections into $\Xast$.
\item $B$ is $w^\ast$-thick.
\end{statements}
\end{thm}

\begin{proof} That (a) implies (b) is trivial. To show that (b) implies (c) we make necessary adjustments in the corresponding proof of Theorem \ref{surcar}. First substitute $A$'s with $B$'s. Then define $C$ by $w^\ast$-closure. Note that $C$ is now a non-$w^\ast$-norming set. Define $Y$ to be the $w^\ast$-closure of $span\:C$. Two cases must be considered:
\begin{enumerate}
\item $Y\neq \Xast$
\item $Y=\Xast$
\end{enumerate}

In the first case let $T$ be the embedding of $Y$ into $\Xast$. Then, since $Y$ is $w^\ast$-closed, $Y=(X/M)^\ast$, where $M$ is the annihilator of $Y$ in $X$. Moreover, $T$ is the adjoint of the quotient map $q:X\rightarrow X/M$.

In the second case $C$ separates points on $X$ and the set $C$ can be used to define a new norm $\|\cdot\|_C$ on $X$ by the formula
\[\|x\|_C=\sup_{c\in C} \|c(x)\|.\]
Then, by the definition of $C$, $\|\cdot\|_C$ is strictly weaker than the original norm. Let $E$ be the completion of $X$ in this weaker norm, let $j$ be the embedding of $X$ into $E$. The adjoint $T$ of $j$ is then injective (since $j$ is dense). Moreover, $T$ is by definition onto $C$ and hence onto $B$. Thus (b) implies (c). 

To show that (c) implies (a), mimicking the corresponding proof of Theorem \ref{surcar} gives the existence of a natural number $j$ such that
\[\overline{(TB_Y)}^{w^{\ast}}\supset\frac{\delta}{j} B_{\Xast}.\]
Now we use that $T$ is an adjoint operator. This gives us that the set $TB_Y$ is $w^\ast$-compact, and hence
\[TB_Y\supset\frac{\delta}{j} B_{\Xast}, \]
which concludes the proof.
\end{proof}

The following corollary is known from \cite{Fs}.

\begin{cor} Suppose $B$ is a subset of a dual space $\Xast$ such that \\ $\overline{span}^{w^{\ast}}(B)=\Xast$. Then the following statements are equivalent.
\begin{statements}
\item There exists a Banach space $Y$ and an injection $T:X\rightarrow Y$ such that $\Tast$ is injective and $\Tast\Yast\supset B$, but $T$ is not invertible.
\item $B$ is $w^\ast$-thin.
\end{statements}
\end{cor}

\begin{cor}\label{wequivcor} In the Banach space setting, the $w^\ast$- surjectivity property and the $w^\ast$- boundedness property are both equivalent to $w^\ast$- thickness. 
\end{cor}

By the Fernandez-Hui-Shapiro theorems (see Theorem \ref{fernandez} and the comments after it), the set of Blaschke-products is a $w^\ast$- thick subset of $H^\infty$. Thus, if $T$ is an adjoint operator from a dual Banach space into $H^\infty$ which is onto the set of Blaschke products, then $T$ is onto $X$. We state this in a way that looks more interesting for applications:

\begin{thm} Let $X$ be an arbitrary Banach space. Let $B$ denote the set of Blaschke-products in $H^\infty$. Suppose $S\in\calL (L_{1}/H_{0}^{1}, X)$ is such that $\Sast(\Xast)\supset B$. Then 
$X$ contains $L_{1}/H_{0}^{1}$ as a closed subspace.
\end{thm}

\section{The Seever property and the Nikodym property}

In \cite[Example 5 p. 18]{DU} an example is given to show that the Nikodym-Grothendieck Boundedness theorem may fail when the measures are not defined on a $\sigma$- algebra but just an algebra.

In \cite{sch} five properties for algebras $\calA$ of sets are discussed. They are as follows:
\begin{ekvivalens}
\item $\calA$ has the Vitali-Hahn-Saks property (VHS) if the Vitali-Hahn-Saks theorem holds on $\calA$.
\item $\calA$ has the Nikodym property (N) if the Nikodym boundedness theorem holds on $\calA$.
\item $\calA$ has the Orlicz-Pettis property (OP) if for every Banach space weak countable additivity implies countable additivity.
\item $\calA$ has the Grothendieck property (G) if $B(\calA)$ is a Grothendieck space.
\item $\calA$ has the Rosenthal property (R) if $B(\calA)$ is a Rosenthal space.
\end{ekvivalens}

It is shown in different papers (see \cite{sch} and \cite{DU2} for references) that (VHS) $\Leftrightarrow$ (N) and (G), that (G) $\Rightarrow$ (OP), that (R) $\Rightarrow$ (G) and that no other implications hold. For some time it was open whether (G) might imply (N). A counterexample was given by M. Talagrand in \cite{tala}. 

Let us say that an algebra has the Seever property (S) if Seever's theorem works on $\calA$. 

\begin{thm} For the following statements about an algebra $\calA$ we have that (a) $\Leftrightarrow$ (b) $\Leftrightarrow$ (c)
\begin{statements}
\item $\calA$ has the Nikodym property (N)
\item $\calA$ has the Seever property (S)
\item The set $\{\chi_A\::A\in\calA\}$ is thick
\end{statements}
\end{thm}  

\section{Some results on thickness in $\calL(X,Y)^\ast$}

Let $X$ and $Y$ be Banach spaces. A very useful set in $\calL(X,Y)^\ast$ is the tensor product $\Xastast\otimes\Yast$. Recall that the action of a functional $\xastast\otimes\yast$ on an operator $T\in\calL(X,Y)$ is defined by $\xastast\otimes\yast(T)=\xastast(\Tast\yast)$. In \cite[Lemma 1.7b p. 268]{hww} it is shown that $ext B_{\Xastast}\otimes ext B_{\Yast}$ is 1-norming for $L(X,Y)$. It need not be $w^\ast$- thick. But often it is.

\begin{lem}\label{tensor} Suppose $A$ and $B$ are $w^\ast $- thick subsets of $\Xastast$ and $\Yast$ respectively. Then $A\otimes B$ is a $w^\ast $- thick subset of $\calL(X,Y)^\ast$.
\end{lem}

\begin{proof} We will use Theorem \ref{wunicar}. Let $(T_n)$ be a sequence in $L(X,Y)$ such that 
\[\sup_n |\xastast\otimes\yast(T_n)|=\sup_n |\xastast(\Tast_n \yast)|<\infty\]
for all $\xastast\in A$ and all $\yast\in B$. Since $A$ is $w^\ast $- thick we conclude that
\[\sup_n \|(\Tast_n \yast)\|<\infty\]
for all $\yast\in B$. Thus $\sup_n |(\Tast_n \yast)(x)|=\sup_n |\yast(T_n x)|<\infty$ for all $\yast\in B$ and every $x\in B_X$. Since $B$ is $w^\ast $- thick, 
\[\sup_n \|T_n x\|<\infty\]
for all $x\in B_X$ and the result follows since $B_X$ is thick.  
\end{proof}

Since by definition no countable set can be $w^\ast$-thick, the extreme points of $B_{l_{1}}$ is a $w^\ast$-thin set. This is in fact a special case of a rather difficult theorem discovered by V.P. Fonf. Recall that a James boundary $J$ for $X$ is a subset of $\Xast$ such that every $x\in X$ attains its norm on $J$. As an example, the set of extreme points of the dual unit ball is a James boundary for any Banach space $X$. 

\begin{thm}\label{fonfc0} If a Banach space $X$ admits a $w^\ast$-thin James boundary $J$, then $X$ contains a copy of $c_0$. 
\end{thm}

For information I will make a list to show how the theorem can be proved with help of different papers. This is the simplest way that I know. 

\begin{proof}
\begin{statements}
\item Note that the restriction to a subspace $Y$ of a James boundary is a James boundary.
\item Put $J=\cup_n A_n$. By Simons' generalization of the Rainwater lemma, there is a sequence $(x_n)$ on $S_X$ which converges weakly to $0$. By the Bessaga-Pe{\l}zcynski selection principle $(x_n)$ can be assumed to be a basic sequence. Let $Y=[x_n]$. We look for $c_0$ inside $Y$. 
\item Let $T$ be the natural embedding of $Y$ in $X$. Put $B_n =\Tast(A_n)$. Then we can show that $J'=\cup B_n$ is a James boundary for $Y$.
\item Show that each $B_n$ is relatively norm-compact as done on page 489 in \cite{Fwe}. Thus $Y$ has a $\sigma$-compact James boundary $J'$.
\item Use Lemma 27 in \cite{fz} to renorm $Y$ equivalently to have a countable James boundary $J''$.
\item Follow the proof of \cite[Theorem 23]{fz} to construct a copy of $c_0$ inside a once more equivalently renormed version of $Y$. This copy is also a copy in $X$. 
\end{statements}
\end{proof}

An interesting result follows from Lemma \ref{tensor} and Theorem \ref{fonfc0}:

\begin{cor} Suppose $X^\ast$ and $Y$ does not contain a copy of $c_0$. Then the set $E=ext\:B_{X^{\ast\ast}}\otimes ext\:B_{\Yast}$ is $w^\ast$-thick in $L(X,Y)^\ast$. 
\end{cor}

\begin{proof} Since $X^\ast$ and $Y$ does not contain a copy of $c_0$ the sets $B_{X^{\ast\ast}}$ and $B_{\Yast}$ are both $w^\ast$-thick. Hence, by Lemma \ref{tensor}, $E$ is $w^\ast$-thick.
\end{proof}

\begin{rem} Note that the set $E$ is not necessarily a James boundary for $L(X,Y)$. But being identical to the set $ext\:B_{K(X,Y)^{\ast}}$ it is a James boundary for $K(X,Y)$.
\end{rem}

By combining the main result from \cite{Fs} with the knowledge of the exposed points of the dual unit ball of $K(X,Y)$ (see e.g. \cite[Theorem 5.1]{ruess}), we obtain the following theorem on $w^\ast$-thickness of $exp\:B_{K(X,Y)^\ast}$.

\begin{thm} Suppose $\Xast$ and $Y$ are separable and $Y$ does not contain a copy of $c_0$. Then 
$exp\:B_{K(X,Y)^\ast}$ is $w^\ast$-thick.
\end{thm}

\begin{proof} Since $\Xast$ is a separable dual it has the RNP and thus doesn't contain a copy of $c_0$. By the main result from \cite{Fs} the sets $A=exp\:B_{X^{\ast\ast}}$ and $B=exp\:B_{Y^\ast}$ are both $w^\ast$-thick. Hence, by Lemma \ref{tensor} $A\otimes B$ is $w^\ast$-thick. But by \cite{ruess} $A\otimes B$ is exactly the set of exposed points of $B_{K(X,Y)^\ast}$.
\end{proof}

\begin{rem} When $\Xast$ and $Y$ both are separable we obtain that $K(X,Y)$ is separable. The point is that $K(X,Y)$ may very well contain $c_0$ even though $\Xast$ and $Y$ doesn't. For example the space $K(l_2)$ contains a copy of $c_0$.
\end{rem}

\begin{cor} Suppose $X$ and $Y$ are separable, reflexive spaces. Then \\ $w^\ast-exp\:B_{K(X,Y)^\ast}$ is $w^\ast$-thick.
\end{cor}

\begin{proof} Use that $w^\ast-exp\:B_{K(X,Y)^\ast}=w^\ast-exp\:B_X \otimes w^\ast-exp\:B_{\Yast}$.
\end{proof}

\begin{cor} Suppose $X$ and $Y$ are separable, reflexive spaces. Then every James boundary of $K(X,Y)$ is $w^\ast$-thick.
\end{cor}

\section{Some questions and remarks}

Suppose a Banach space contains a thin, fundamental set. Then, by definition of such a set there exists a $w^\ast$- null sequence $(\xast_n)$ on $S_{\Xast}$. Thus, the Josefsson-Nissenzweig theorem is just a triviality in such spaces. Of course, every separable Banach space contains a thin, even norming, set (take any dense countable subset of $B_X$). Also it is immediate that any Banach space containing a complemented separable subspace contains a thin, norming set. Thus, WCG spaces contains thin, norming sets.

\begin{ques} Does every Banach space have a bounded, fundamental, thin set? Does every Banach space have a bounded, norming thin set?
\end{ques}

Motivated by Theorem \ref{surcar} we ask the following question:

\begin{ques} Suppose $A\subset B_X$ is thin. Is there a Banach space $Y$ and a semi-embedding $T:Y\rightarrow X$ such that $TY\supset A$ but $T$ is not onto.
\end{ques}

Fernandez, Hui and Shapiro has asked (in our language) whether the Blaschke-products is a thick subset in $H^\infty$ (not only $w^\ast$- thick). We formulate an extended question:

\begin{ques} Is the set of inner functions (Blaschke-products) a thick subset in $H^\infty$? Is the set of interpolating Blaschke-products thick or $w^\ast$- thick? 
\end{ques}

By a theorem of Mooney (see \cite[p. 206-207]{garnett}), the pre dual of $H^\infty$ is weak sequentially complete. Thus, it doesn't contain a copy of $c_0$. Hence, any James boundary in $H^\infty$ is $w^\ast$-thick. In light of V.P Fonf's theorem (Theorem \ref{fonfc0}) it is natural to ask:

\begin{ques} Is the set of inner functions (Blaschke-products, interpolating Blaschke-products) a James boundary in $H^\infty$?
\end{ques}

If so, Fernandez' and Shapiro's results would follow as special cases of Theorem \ref{fonfc0}. 

We end this paper by giving a list of sets for which results on thickness are known:
\begin{thm}
\begin{statements}
\item Any James-boundary of a Banach space not containing $c_0$ is $w^*$-thick (Fonf \cite{Fwe}). Thus, in particular $ext\:B_X$ is thick for every reflexive Banach space $X$.
\item If $X$ is a separable Banach space not containing $c_0$, then $exp\:B_{\Xast}$ is $w^*$-thick (Fonf \cite{Fs}).
\item The set of characteristic functions in $B(\calA)$ when $\calA$ has the Nikodym (Seever) property is thick (Nikodym-Grothendieck \cite{DU}). Thus $ext\:B_{l_\infty}$ is thick.
\item The set of Inner functions and Blaschke products in $H^\infty$ are $w^\ast$-thick and norming for the dual (\cite{fernandez}, \cite{hui}). 
\item The tensor product of two $w^\ast$-thick sets in $\Xastast$ and $\Yast$ is a $w^\ast$-thick subset in $L(X,Y)^\ast$ (this paper).
\item Suppose $X$ and $Y$ are separable, reflexive spaces. Then every James boundary of $K(X,Y)$ is $w^\ast$-thick.
\end{statements}
\end{thm}


\end{document}